\documentclass[12pt]{article}
\usepackage{amsmath, latexsym}
\title{Counting Permutations by Their Runs Up and Down}
\author{E. Rodney Canfield \thanks{Research supported by the
                             NSA Mathematical Sciences Program}\\
  University of Georgia\\  Athens, GA 30602\\ {\small \tt erc@cs.uga.edu}
  \and Herbert S. Wilf\\University of Pennsylvania\\Philadelphia,
PA 19104-6395\\
   {\small \texttt{wilf@math.upenn.edu}}}
\newtheorem{theorem}{Theorem}
\newtheorem{lemma}{Lemma}
\newcommand{\eqdef}{\, =\kern -12.7pt\raise 6pt\hbox{{\tiny\textrm{def}}}\,\,}
\begin{document}
\maketitle

\begin{abstract}
We find a formula for the number of permutations of $[n]$ that have
exactly $s$ runs up and down. The formula is at once terminating,
asymptotic, and exact.
\end{abstract}

\section{Introduction}
We will say that a \textit{run} of a permutation $\sigma$ is a
maximal interval of consecutive arguments of $\sigma$ on which the
values of $\sigma$ are monotonic. If the values of $\sigma$ increase
on the interval then we speak of a run up, else a run down.
Throughout this paper we will use the unqualified term \textit{run}
to mean either a run up or a run down. These runs have been called
\textit{sequences} by some other authors, and have been called
\textit{alternating runs} by others. For example, the permutation
\[(723851469)\] has four runs, viz. 72, 238, 851, 1469. We let
$P(n,s)$ denote the number of permutations of $n$ letters that have
exactly $s$ runs. Here are the first few values of $P(n,s)$:
$$
\begin{matrix}
   n\backslash k & 1 &  2 &  3 &  4 \cr
        2        & 2                \cr
        3        & 2 &  4           \cr
        4        & 2 & 12 & 10      \cr
        5        & 2 & 28 & 58 & 32
\end{matrix}
$$
There is a large literature devoted to this $P(n,s)$, which we will
survey briefly in section \ref{sec:surv}. But although a number of
recurrences and generating functions etc. are known, it does not
seem to have been noticed that an interesting exact formula of the
kind we present in this paper exists. Carlitz \cite{ca2} has derived
an exact formula for $P(n,s)$, but that one is not at the same time
an asymptotic formula.  The Carlitz formula is discussed further in
the final section.

Our approach to this problem differs from  previous studies in that
we concentrate on the generating functions $u_s(x)$, defined for
each fixed $s\ge 1$ by \[u_s(x)=\sum_{n\ge 2}P(n,s)x^n,\] whereas
most earlier work has dealt with generating functions for fixed $n$.
By finding the form of these generating functions we will be able
to exhibit a formula for $P(n,s)$ which is simultaneously
\begin{itemize} \item exact, and \item terminating, and \item
asymptotic, for fixed $s$ and $n\to\infty$.
\end{itemize}
To our knowledge, the asymptotic behavior of the $P(n,s)$ has not
been previously explored.

The formula that we will find is of the form \begin{equation}
P(n,s)=\frac{s^n}{2^{s-2}}-\frac{(s-1)^n}{2^{s-4}}+\psi_2(n,s)(s-2)^n+\dots+
\psi_{s-1}(n,s),\qquad (n\ge 2)
\end{equation}
in which each $\psi_i(n,s)$ is a polynomial in $n$ whose degree in
$n$ is $\lfloor{i/2}\rfloor$.
\section{Outline of this paper}
In section \ref{sec:us} we will find the generating functions
$u_s(x)=\sum_nP(n,s)x^n$, as a rational function. Since the
denominator will appear in completely factored form, we can write
out, in section \ref{sec:pform}, a formula for $P(n,s)$ of the type
described above.

Interestingly, the formula will be, in that section, uniquely
determined except for the coefficient of the leading term! That is,
we will show in that section, that for fixed $s$ we have $P(n,s)=
K(s)s^n+\dots$, but $K(s)$ will be, for the moment, unknown.

In section \ref{sec:ks} we begin the task of determining the
multiplicative factor $K(s)$. Surprisingly, although the tools that
will have been used up to that point will be entirely analytical in
nature, the determination of $K(s)$ will be done by an
``almost-bijection.'' We will show that $P(n,s)$ is, for fixed $s$,
asymptotic to the number of $s$-tuples of pairwise-disjoint subsets
of $[n]$, each of cardinality $\ge 2$, and the asymptotic behavior
of the latter is easily found.

The combination of the former analytical results and the latter
bijective argument results in the complete formula for $P(n,s)$.

\section{Finding the $u_s(x)$ functions} \label{sec:us} The recurrence
formula for the numbers $P(n,s)$ is well known and is due to Andr\'e
\cite{an}, \begin{equation} \label{eq:andre}
P(n,s)=sP(n-1,s)+2P(n-1,s-1)+(n-s)P(n-1,s-2),\qquad (n\ge 3) \end{equation}
with $P(2,s)=2\delta_{s,1}$. {From} this
recurrence one finds easily a recurrence for the generating
functions $u_s(x)\eqdef\sum_nP(n,s)x^n$, viz.
\begin{equation}
\label{eq:ugf}
(1-sx)u_s(x)=2xu_{s-1}(x)+x^2u_{s-2}'(x)-(s-1)xu_{s-2}(x),\qquad(s\ge 2)
\end{equation}
with $u_1(x)=2x^2/(1-x),u_0(x)=0$. The next three of these functions
are \begin{eqnarray*} u_2(x)&=&\frac{4x^3}{(1-x)(1-2x)},\\
u_3(x)&=&\frac{2x^4(5-6x)}{(1-3x)(1-2x)(1-x)^2},\\
u_4(x)&=&\frac{4x^5(8-29x+24x^2)}{(1-4x)(1-3x)(1-2x)^2(1-x)^2}.
\end{eqnarray*}

We will find the general form of these functions, and from that will
follow the desired formulas for $P(n,s)$.
\begin{theorem}
\label{thm:uform} We have, for each $s=1,2,3,\dots,$
\begin{equation}
\label{eq:uform}
u_s(x)=\frac{\Phi_s(x)}{(1-sx)(1-(s-1)x)(1-(s-2)x)^2(1-(s-3)x)^2\cdots
(1-x)^{\lfloor{(s+1)/2\rfloor}}},\end{equation} where $\Phi_s(x)$
is
a polynomial of degree
$1+\left\lceil{\frac{s(s+2)}{4}}\right\rceil$. The degree of the
denominator is $\left\lceil{\frac{s(s+2)}{4}}\right\rceil$, which
is
exactly 1 less than the degree of the numerator, for all $s\ge 1$.
\end{theorem}
\subsection{Proof of Theorem \ref{thm:uform}}

The proof is by a straightforward, though tedious, substitution of
the form (\ref{eq:uform}) into the recurrence (\ref{eq:ugf}) to find
a recurrence for the numerator polynomials $\Phi_s(x)$. This will
establish that they are indeed polynomials and will provide the
claimed degree estimates. We will do this by putting every term over
the common denominator \begin{eqnarray*}
\Delta_s(x)&=&(1-sx)(1-(s-1)x)(1-(s-2)x)^2(1-(s-3)x)^2\cdots
(1-x)^{\lfloor{(s+1)/2\rfloor}}\\
 &\eqdef&
\prod_{i=0}^{s-1}(1-(s-i)x)^{\epsilon_i},
\end{eqnarray*}
where we have written $\{\epsilon_i\}_{i\ge
0}=\{1,1,2,2,3,3,4,4,\dots\}$.

For technical reasons it will be useful to rewrite the recurrence
(\ref{eq:ugf}) in the form
\begin{eqnarray}
u_s(x)&=&\frac{2xu_{s-1}(x)}{(1-sx)}+\frac{x^2u_{s-2}'(x)}{(1-sx)}-\frac{(s-1)xu_{s-2}(x)}{(1-sx)}\
\nonumber\\
&=&\frac{2x\Phi_{s-1}(x)}{(1-sx)\Delta_{s-1}(x)}+\frac{x^2\Phi_{s-2}'(x)}{(1-sx)\Delta_{s-2}}-
\frac{x^2\Phi_{s-2}}{(1-sx)\Delta_{s-2}}\frac{\Delta_{s-2}'(x)}{\Delta_{s-2}(x)}\nonumber\\
&&\qquad\qquad
-\frac{(s-1)x\Phi_{s-2}(x)}{(1-sx)\Delta_{s-2}(x)}\nonumber\\
&=&\frac{1}{\Delta_s(x)}\bigg\{\frac{2x\Phi_{s-1}(x)\Delta_s(x)}{(1-sx)\Delta_{s-1}(x)}
+\frac{x^2\Phi_{s-2}'(x)\Delta_s(x)}{(1-sx)\Delta_{s-2}}-
\frac{x^2\Phi_{s-2}(x)\Delta_s(x)}{(1-sx)\Delta_{s-2}}\frac{\Delta_{s-2}'(x)}{\Delta_{s-2}(x)}\nonumber\\
&&\qquad\qquad
-\frac{(s-1)x\Phi_{s-2}(x)\Delta_s(x)}{(1-sx)\Delta_{s-2}(x)}\bigg\}\label{eq:urcr}
\end{eqnarray}
Each of the four terms inside the braces is a polynomial in $x$
whose degree we will now find.

Consider the ratio
\begin{eqnarray*}
\frac{\Delta_s(x)}{(1-sx)\Delta_{s-1}(x)}
&=&\frac{\prod_{j=0}^{s-1}(1-(s-j)x)^{\epsilon_j}}{(1-sx)\prod_{j=0}^{s-2}(1-(s-1-j)x)^{\epsilon_j}}\\
&=&\frac{\prod_{j=0}^{s-1}(1-(s-j)x)^{\epsilon_j}}{(1-sx)\prod_{j=1}^{s-1}(1-(s-j)x)^{\epsilon_{j-1}}}\\
&=&\prod_{j=1}^{s-1}(1-(s-j)x)^{\epsilon_j-\epsilon_{j-1}}=\prod_{j\,\mathrm{even};\,2\le
j\le s-1}(1-(s-j)x),
\end{eqnarray*}
which is a polynomial of degree $\lfloor{(s-1)/2}\rfloor$.

 It follows that
\begin{eqnarray*}
\frac{\Delta_s(x)}{(1-sx)\Delta_{s-2}(x)}&=&\left(\frac{\Delta_s(x)}{(1-sx)\Delta_{s-1}(x)}\right)
\left(\frac{\Delta_{s-1}(x)}{(1-(s-1)x)\Delta_{s-2}(x)}\right)(1-(s-1)x)\\
&=&\prod_{j\,\mathrm{even};\,2\le j\le
s-1}(1-(s-j)x)\prod_{j\,\mathrm{even};\,0\le j\le s-2}(1-(s-1-j)x)\\
&=&\prod_{j\,\mathrm{even};\,2\le j\le
s-1}(1-(s-j)x)\prod_{j\,\mathrm{odd};\,1\le j\le s-1}(1-(s-j)x)\\
&=&\prod_{j=1}^{s-1}(1-(s-j)x), \end{eqnarray*} is a polynomial in
$x$ of degree $s-1$.

We can now deal with the third term of the four inside the braces
in
(\ref{eq:urcr}). Since
\[\frac{\Delta_{s-2}'(x)}{\Delta_{s-2}(x)}=-\sum_{j=2}^{s-1}\frac{\epsilon_{j-2}(s-j)}{1-(s-j)x},\]
we have
\begin{eqnarray*}
\frac{x^2\Phi_{s-2}(x)\Delta_s(x)}{(1-sx)\Delta_{s-2}}\frac{\Delta_{s-2}'(x)}{\Delta_{s-2}(x)}&=&
x^2\Phi_{s-2}(x)\left(\prod_{j=1}^{s-1}(1-(s-j)x)\right)\left(\sum_{j=2}^{s-1}
\frac{-\epsilon_{j-2}(s-j)}{1-(s-j)x}\right). \end{eqnarray*} If
$d(s)$ denotes the degree of $\Phi(s)$, then this last member is a
polynomial in $x$ of degree $2+d(s-2)+s-2=d(s-2)+s$.

We have now shown that each of the four terms inside the braces in
(\ref{eq:urcr}) is a polynomial in $x$. Their respective degrees are
\[d(s-1)+1+\lfloor{(s-1)/2}\rfloor,d(s-2)+s,d(s-2)+s,d(s-2)+s.\]
Hence we have
$d(s)=\max{(d(s-1)+\lfloor{(s+1)/2}\rfloor,d(s-2)+s)}$, with
$d(2)=3$ and $d(3)=5$.

It is remarkable that this difference equation has a simple
solution. Its solution is
\[d(s)=1+\left\lceil{\frac{s(s+2)}{4}}\right\rceil,\] as can easily
be checked, and in fact all four terms inside the braces in
(\ref{eq:urcr}) have the same degree! This completes the proof of
the Theorem.
\section{The formula for $P(n,s)$}\label{sec:pform}
{From} the partial fraction expansion of (\ref{eq:uform}) we find
at
once that
\begin{equation}
\label{eq:pform}
P(n,s)=\psi_0(n,s)s^n+\psi_1(n,s)(s-1)^n+\psi_2(n,s)(s-2)^n+\dots+\psi_{s-1}(n,s),\qquad
(n\ge 2)
\end{equation}
where each $\psi_i(n,s)$ is a polynomial in $n$ of degree at most
$\lfloor{i/2}\rfloor$, and it remains to find these polynomials. To
do this we substitute (\ref{eq:pform}) into the recurrence
(\ref{eq:andre}) and match the coefficients of each term $(s-i)^n$.
The result of this substitution is that
\[(s-i)\psi_i(n,s)=s\psi_i(n-1,s)+2\psi_{i-1}(n-1,s-1)+(n-s)\psi_{i-2}(n-1,s-2).\]
Perhaps the best way to find these $\psi$'s explicitly is to assume
a solution in the form of a polynomial in $n$ of degree
$\lfloor{i/2}\rfloor$ and solve for the coefficients of that
polynomial. We can begin with $\psi_{-1}(n,s)=0$ and
$\psi_0(n,s)=K(s)$ (since $\psi_0$ is of degree zero in $n$) where
$K$ is to be determined. We then find that
\[\psi_1(n,s)=-2K(s-1),\qquad \psi_2(n,s)=\frac14 K(s-2)(s+8-2n)\]
\[\psi_3(n,s)=\frac12 K(s-3)(2n-s-3),\quad \psi_4(n,s)=\frac{1}{32}
K(s-4)(4n^2-4n(s+8)+s^2+15s+32)\] For example, for $s=4$ we find $$
u_4(x) =
  \frac{1/4}{1-4x}
+ \frac{(-1)}{1-3x} + \frac{(-1/2)}{(1-2x)^2} + \frac{7/2}{1-2x} +
\frac{2}{(1-x)^2} + \frac{(-9)}{1-x} + 2x + 19/4.
$$
{From} this it follows
$$
P(n,4) = 4^{n-1} - 3^n + (6-n)2^{n-1} + (2n-7), ~~ (n \ge 2).
$$

\section{The factor $K(s)$}\label{sec:ks} We have now described the
formula for $P(n,s)$ completely except for the multiplicative factor
$K(s)$. It remains to show that $K(s)=2^{-(s-2)}$.  For this, it would
suffice to prove the next Theorem for fixed $s$ and $n\rightarrow\infty$;
since the proof is applicable to a larger range of $s$, we state it
in that manner:
\begin{theorem}
\label{thm:K} Let $\epsilon>0$, and $\{(n,s)\}$ be an infinite
sequence of pairs such that $n \rightarrow \infty$ and $s \le
(1+\epsilon)^{-1}n/\log n$.  Then, \begin{equation} \label{eq:asym}
P(n,s)\sim \frac{1}{2^{s-2}}s^n.
\end{equation}
\end{theorem}

\subsection{Proof of Theorem \ref{thm:K}}

To fix ideas, we will do this by showing that the number
$\hat{P}(n,s)$ of permutations of $n$ letters, with $s$ runs,
\textit{the first of which is a run up}, is $\sim s^n/2^{s-1}$.
Evidently the number for which the first run is down will be the
same, and the desired result will follow. Henceforth we will always
assume that the first run is a run up.  There are two steps to the
proof.  In the first step, we show that the set of permutations
counted by $\hat{P}(n,s)$ can be put into bijection with $s$-tuples
of subsets $(S_1,\dots,S_s)$  (each $S_i\subseteq [n]$) satisfying
certain properties.  In the second part of the proof, we introduce
a
function called $\Phi$ whose domain is the Cartesian product of
these $s$-tuples with a set of cardinality $2^{s-1}$, and whose
range is a set of size $s^n$.  We prove that this function $\Phi$
is
an \textit{injection}. Although we have no succinct description of
the image of this injection, we are able to show that for $(n,s)$
in
the range hypothesized by the theorem the image is asymptotically
all of the range set.

\subsection{First part of the proof}
Let $\Pi(n,s)$ be the set of all $n$-permutations with $s$ runs up
and down, the first of which is up. Let $\tilde{\Pi}(n,s)$ be the
collection of all $s$-tuples $(S_1,\dots,S_s)$ of nonempty subsets
of $[n]$ which are \textit{almost} pairwise disjoint, in that
\begin{equation}\label{eq:disj} |S_i\cap S_j|=
\begin{cases}
1,&\mathrm{if}\ j=i+1\ \mathrm{and}\ 1\le i< s;\\ 0,&\mathrm{else}
\end{cases} \end{equation} Further we require that
\begin{equation}\label{eq:big} |S_i|\ge 2,\ \forall i,
\end{equation}
and that \begin{eqnarray} \max{(S_i)}&=&\max{(S_{i+1})}\in S_i\cap
S_{i+1}\qquad (\forall\ \mathrm{odd}\ i)\nonumber\\
\min{(S_i)}&=&\min{(S_{i+1})}\in S_i\cap S_{i+1}\qquad (\forall\
\mathrm{even}\ i).\label{eq:cont} \end{eqnarray} \begin{lemma} The
number of $s$-tuples of subsets of $[n]$ that satisfy
(\ref{eq:disj})--(\ref{eq:cont}) is equal to the number of
permutations  of $[n]$ with $s$ runs, the first of which is up.
\end{lemma}
Indeed to reconstruct the permutation from the $s$-tuple of sets,
we
first sort each of the sets, the first in increasing order, the
second decreasing, etc., then merge the sets, and finally delete one
element of each of the adjacent duplicates that appear. $\Box$

Hence it suffices to show that the number of $s$-tuples of subsets
that satisfy (\ref{eq:disj})--(\ref{eq:cont}) is $\sim s^n/2^{s-1}$.

\subsection{Defining the function $\Phi$} By a \textit{choice sequence}
$\mathbf{h}=(h_1,\dots,h_{s-1})$ we mean an $s-1$-tuple where each
$h_i$ is either equal to $i$ or to $i+1$. The set of all such choice
sequences will be $H_s$.  The function to be constructed is a mapping
\[\Phi: H_s\times \tilde{\Pi}(n,s)\rightarrow \{(T_1,T_2,\dots,T_s):
\forall i,T_i\subseteq [n]\}.\]

Let $\mathbf{h}\in H_s$, and let $(S_1,\dots,S_s)$ be a family of
subsets satisfying (\ref{eq:disj})--(\ref{eq:cont}). For
each $i=1,\dots,s-1$, let $e_i$ be the unique element that belongs
to $S_i\cap S_{i+1}$. These $e_i$'s are all different, since
$e_i=e_j$ with $i<j$ would imply that $S_i\cap S_{j+1}$ is nonempty,
contradicting (\ref{eq:disj}). Perform the following $s-1$ delete
operations: for each $i=1,\dots,s-1$, delete the element $e_i$ from
the set $S_{h_i}$.  The resulting $s$-tuple of sets remaining after
these deletions is, by definition,
$\Phi\bigl(\mathbf{h},(S_1,\dots,S_s)\bigr)$.

The image of this mapping does not include all $s$-tuples of sets,
as the following Lemma shows.
\begin{lemma} If $(T_1,\dots,T_s)$ is in the image of $\Phi$ then
\begin{enumerate} \item the $T_i$'s are pairwise disjoint, and \item
the union of the $T_i$'s is $[n]$.
\end{enumerate}
$\Box$
\end{lemma} It is possible for some of the $T_i$'s to be empty. We
remark that the number of $s$-tuples
$(T_1,\dots,T_s)$ in which the $T_i$'s are pairwise disjoint and
whose union is $[n]$ is $s^n$.

\subsection{The mapping $\Phi$ is injective} The way we prove this
assertion is to give a \textit{reconstruction} algorithm.  The algorithm
begins with an $s$-tuple $(T_1,\dots,T_s)$ of subsets which putatively
belongs to the image of $\Phi$.  It attempts to reconstruct the preimage.
It will be clear from the algorithm that the preimage can be only
one thing, if it exists at all.  There is one ``early exit''
point in the algorithm where the search for a preimage is abandoned,
because it obviously does not exist.  If the algorithm executes all
the way to finish, then it will have found the only possible 
candidate for a preimage.  However, it is still possible that the
$s$-tuple of sets found at the end will not satisfy one of the
required conditions (\ref{eq:disj})--(\ref{eq:cont}).

\begin{lemma} The mapping $\Phi$ is injective.\end{lemma} \noindent{\textbf
Proof.} Let $(T_1,\dots,T_s)$ be an $s$-tuple of pairwise disjoint
(possibly empty) sets whose union is $[n]$.  Here is the reconstruction
algorithm:
\begin{enumerate}
\item (Find consecutive unions) It is easy to see that if the
$T$'s are in the image of $\Phi$, then \[T_i\cup
T_{i+1}=S_i\cup S_{i+1},\quad 1\le i<s.\]  So if one of the inequalities
\[|T_i\cap T_{i+1}|\ge 3\quad 1\le i<s,\] fails, then the reconstruction
fails and no preimage exists. Otherwise we can reconstruct all of
the unions $S_i\cup S_{i+1}$.
\item (Reconstruct the set of deleted elements) Put $e_1=\max{(S_1\cup
S_2)},e_2=\min{(S_2\cup S_3)},\dots$.
\item (Recover the choice sequence $\mathbf{h}$) For each $i=1,\dots,s-1$,
since $e_i\in T_i\cup T_{i+1}$, and because the $T_i$'s are pairwise
disjoint, there will be exactly one index, $h_i$, say, such that $h_i\in\{i,i+1\}$
and $e_i\notin T_{h_i}$.
\item (Re-insert the elements that were deleted) For each $i$, $1\le
i<s$, insert the element $e_i$ into the set $T_{h_i}$.
\end{enumerate}
If the reconstructed sets $(S_1,\dots,S_s)$ satisfy
(\ref{eq:disj})--(\ref{eq:cont}) then we have found the unique
preimage. Otherwise no preimage exists. $\Box$

Thus if $\hat{P}(n,s)$ is the number of permutations of $n$ letters
with $s$ runs, the first of which is up, then we have shown that
\begin{equation} \label{eq:upper} 2^{s-1}\hat{P}(n,s)\le s^n.
\end{equation}

\subsection{When does the algorithm terminate without a preimage?}
If the reconstruction algorithm does not early exit in step 1,
yet fails to find a preimage, then one of the conditions
(\ref{eq:disj})--(\ref{eq:cont}) is not satisfied. We will now visit
each of these in turn to see when it might fail \begin{enumerate}
\item (Can (\ref{eq:disj}) fail?) The intersections $T_i\cap
T_{i+1}$ were all empty before the insertions; however, the
operation, ``insert element $e_i$ into $T_{h_i}$" either added an
element of $T_i$ to the set $T_{i+1}$, or vice-versa. That operation
alone caused the two adjacent sets to have intersection 1. The only
other insertion which could have affected $T_i$ is the one which
involves element $e_{i-1}$. If that operation increased the size of
$T_i$, then it did so by inserting an element from $T_{i-1}$, which
element could not possibly be present in $T_{i+1}$. Thus, the only
other insertion which could possibly affect the set $T_i$ will have
no effect on the cardinality of $T_i\cap T_{i+1}$. Likewise, the
only other operation which can possibly affect the cardinality of
$T_{i+1}$ will have no effect on the cardinality of $T_i\cap
T_{i+1}$. So, the intersection $S_i\cap S_{i+1}$ will always have
size 1, as required.
\item (Will $S_i\cap S_j=\emptyset$ when $j>i+1$?) Only the case $j
= i+2$ is not obvious. If $S_i\cap S_{i+2}$ is not empty, then during
reconstruction some element originally belonging to $T_{i+1}$ was
inserted into both $T_i$ and $T_{i+2}$. (Any element originally in
$T_i$ cannot end up in $S_{i+2}$, and
vice-versa.) This means that some element $e\in T_{i+1}$ is both the
maximum of $T_i\cup T_{i+1}$, as well as the minimum of $T_{i+1}\cup
T_{i+2}$; (or the other way around). But \[\max{(T_{i+1})}\le
\max{(T_i\cup T_{i+1})}=e_i=\min{(T_{i+1}\cup T_{i+2})} \le
\min{(T_{i+1})},\] so $T_{i+1}$ has just one element, $e_i$, and
that element lies between the surrounding sets $T_i$ and $T_{i+2}$.
This is how the reconstructed sets can fail to satisfy
(\ref{eq:disj}).
\item (Can (\ref{eq:big}) fail?) Yes, if one of the sets $T_i$ has
 cardinality 0 or 1, then it is possible that the cardinality of 
the reconstructed $S_i$ may be less than 2.
\item (Can (\ref{eq:cont}) fail?) No. By the nature of the reconstruction,
the $S_i$'s always have this property.
\end{enumerate}

We can now prove
\begin{lemma}
If in the given sequence $T=(T_1,\dots,T_s)$, all sets have
cardinalities at least 2, then $T$ has a preimage under $\Phi$.
\end{lemma}
For then the unions $T_i\cup T_{i+1}$ have size 4 or more, so we
don't terminate the reconstruction at Step 1. The only other two
possible failures --- when an intersection $S_i\cap S_{i+2}$ was
nonempty, or one of the $S_i$ was too small
--- could both be traced back to a set $T_i$ which had size 0 or 1.
$\Box$

A crude lower estimate from Bonferroni's inequalities tells us that
the number of $s$-tuples $T$ that are pairwise disjoint, with union
equal to $[n]$, and with all cardinalities $\ge 2$ is at least
\[s^n-(n+s)(s-1)^{n-1}.\] The reason: $s(s-1)^n$ is an upper bound
on the number of $T$'s for which some component is the empty set;
and $ns(s-1)^{n-1}$ is an upper bound on the number of $T$'s for
which some component has cardinality one; then, $s(s-1)+ns<s(n+s)$.
Hence
\begin{equation}
\label{eq:lower} 2^{s-1}\hat{P}(n,s)\ge s^n-s(n+s)(s-1)^{n-1},
\end{equation} which, taken together with (\ref{eq:upper}) completes
the proof of (\ref{eq:asym}), since our hypothesis on the pairs
$(n,s)$ implies $$ (n+s)(s-1)^{n-1} \le 2n(s-1)^{n-1} = o(s^n).
$$

\vskip 5pt


\section{Survey of the literature} \label{sec:surv}  Andr\'e was the
first to study \cite{an} the runs up and down of permutations, and
the fundamental recurrence
(\ref{eq:andre}) is due to him. His paper includes a table of
$P(n,s)$ through $n=8$, with one error in the final row. 
A great deal of information about $P(n,s)$ is found
in vol. 3 of \cite{kn} (see
particularly ex. 15, 16 of sec. 5.1.3).  Comtet \cite{comtet}
devotes an extended exercise, see page 260, to the topic.  The
two variable generating function given there, however, 
is incorrect.  A correct version appears in the discussion 
accompanying  sequence A059427
of \cite{sloane}.

Carlitz
\cite{ca1,ca2,ca3} visited this subject several times. In \cite{ca1}
he gives the two-variable generating function $$
\sum_{n=2}^{\infty}\frac{z^n}{n!} (1-x^2)^{-n/2} \sum_{s=1}^{n-1}
P(n+1,s) x^{n-s} = \frac {(1-x)((1-x^2)^{1/2}+\sin(z))^2}
{(1+x)(x-\cos(z))^2}, $$ and in \cite{ca2} he finds an explicit
formula for $P(n,s)$ and information about an associated polynomial
sequence. There is something wrong with the final formulas of this
latter work, however; these formulas suggest $P(8,s) =
0,2,250,2516,7060,7562,2770$; whereas, in fact, $P(8,s) = 2, 252,
2766, 9576, 14622, 10332, 2770$.  (Empirically, his formula always
gives the right value for $P(n,n-1)$.) Further evidence that
something is amiss concerns the auxiliary quantity
$\overline{K}_{n,j}$; the summation formula given for this quantity
does not give the values displayed in the table.  Whether the
problem can be easily repaired, we have not investigated.

More recently, B\'{o}na and
Ehrenborg \cite{be} have proven log-concavity: $P(n,s)^2 \ge
P(n,s-1)P(n,s+1)$.  In the later book \cite{bona}, the stronger
assertion, that $P_n(x) \eqdef \sum_s P(n,s)x^s$ has all its roots
real and negative, is made.  A proof of this can be based on the
relation $$ P_n(x) = (x-x^3)P_{n-1}'(x) + \bigl( (n-2)x^2+2x \bigr)
P_{n-1}(x),$$ which itself is a consequence of the basic recursion
(\ref{eq:andre}). This implies, once it is established that the
variance becomes infinite with $n$, that the numbers $P(n,s)$
satisfy a central limit theorem. (That is, are asymptotically
normal.)  Due to log-concavity, one may deduce (see Theorem 4 of
\cite{bender}) a local limit theorem. This leads to an asymptotic
formula for $P(n,s)$ for $s$ in a different range than in our
Theorem \ref{thm:K}.

\end{document}